\documentclass[conference]{IEEEtran}
\IEEEoverridecommandlockouts

\usepackage{pgfplots}
\usepackage{pgfplotstable}
\pgfplotsset{compat=newest}
\usepackage{cite}
\usepackage{amsmath,amssymb,amsfonts}
\usepackage{algorithmic}
\usepackage{graphicx}
\usepackage{textcomp}
\usepackage{xcolor}
\usepackage[font=footnotesize]{caption}
\usepackage{subcaption}
\usepackage{url}  
\usepackage{xurl}  
\usepackage{amsmath}
\usepackage{bm}  

\def\BibTeX{{\rm B\kern-.05em{\sc i\kern-.025em b}\kern-.08em
    T\kern-.1667em\lower.7ex\hbox{E}\kern-.125emX}}
\begin{document}

\title{



GPU-Accelerated Optimization Solver for Unit Commitment in Large-Scale Power Grids 

}

\author{\IEEEauthorblockN{Hussein Sharadga}
\IEEEauthorblockA{~~~~\textit{School of Engineering~~~~} \\
\textit{~~~~~Texas A\&M International University}\\
Laredo, Texas, USA \\
hussein.sharadga@tamiu.edu}
\and
\IEEEauthorblockN{Javad Mohammadi}
\IEEEauthorblockA{~~~~~~~~\textit{Civil, Architectural and Environmental Engineering~~~~~~~~} \\
\textit{The University of Texas at Austin}\\
Austin, Texas, USA \\
javadm@utexas.edu}\\
}
\maketitle

\IEEEpubidadjcol

\begin{abstract}

This work presents a GPU-accelerated solver for the unit commitment (UC) problem in large-scale power grids. The solver uses the Primal–Dual Hybrid Gradient (PDHG) algorithm to efficiently solve the relaxed linear subproblem, achieving faster bound estimation and improved crossover and branch-and-bound convergence compared to conventional CPU-based methods. These improvements significantly reduce the total computation time for the mixed-integer linear UC problem. The proposed approach is validated on large-scale systems, including 4224-, 6049-, and 6717-bus networks with long control horizons and computationally intensive problems, demonstrating substantial speed-ups while maintaining solution quality.

\end{abstract}

\begin{IEEEkeywords}
Power Grid Optimization, Large-Scale Optimization, Unit Commitment, GPU, Optimal Power Flow (OPF), Primal-Dual Hybrid Gradient (PDHG) Algorithm, Primal-Dual Linear Programming (PDLP) 
\end{IEEEkeywords}

\vspace{-.4cm}
\section{Introduction}

\subsection{Motivation}
Although the electrical power grid is widely regarded as one of the landmark achievements of 20th-century engineering \cite{b1}, efficiently managing power distribution continues
to be a significant challenge. With the growing complexity of modern power systems, ensuring system reliability and resilience has become increasingly important. To support these objectives, the U.S. Department of Energy (DOE) promotes programs aimed at modernizing grid operations through the adoption of advanced technologies and methodologies \cite{b2}.

Optimal power flow (OPF) plays a critical role in grid management by determining optimal generator outputs, load allocations, energy storage operation, and other control settings \cite{b3, SCOPF-Javad}. Security-constrained OPF (SCOPF) further strengthens system reliability by explicitly considering potential contingencies within the operational planning process \cite{b4}.

Historically, SCOPF methods were developed under limited computational capabilities, often employing linearized models such as DCOPF that ignore voltage magnitudes and reactive power effects \cite{b5}. Despite improvements in algorithms and computational power, many existing tools continue to rely on these simplifications, making it challenging to solve the full AC power flow problem \cite{b2}. To overcome these limitations, the DOE launched the Grid Optimization (GO) challenges, which focus on solving the full AC power flow problem while accounting for unit commitment and security requirements \cite{b2}.

\subsection{Prior Work on the Grid Optimization (GO) Challenge}
The GO Challenge series focuses on large-scale nonlinear optimization tasks involving unit commitment decisions, which are known to be computationally demanding. To manage this complexity, several research efforts have introduced advanced optimization techniques aimed at improving scalability and solution performance. For example, Chevalier \cite{ref4} proposed a parallel Adam-based optimization method that speeds up both backpropagation and variable projection through parallel computation, allowing more efficient processing of large problem instances. In related work, \cite{ref3} showed that enforcing AC power flow constraints, reserve requirements, and unit commitment  can produce reliable solutions for AC-based unit commitment challenges. 
Additionally, Bienstock and Villagra \cite{refX} introduced a linear relaxation framework for ACOPF based on cutting-plane methods, enabling the computation of robust lower bounds through outer-envelope cuts and refined cut management procedures. 
Holzer et al. \cite{refY} evaluated multiple solvers submitted to GO Challenge 3, which addressed a multi-period unit commitment problem with AC network modeling and topology switching, illustrating the comparative performance of different optimization strategies at scale.

In prior work \cite{TPEC, IEEE-TIA}, we proposed a solution strategy for the full AC power flow challenge, decomposing it into two sequential subproblems: a DC model that handles the binary unit commitment decisions and feasibility constraints, and an AC model that solves continuous variables under all constraints using a nonconvex solver. This decomposition enables efficient handling of computational workloads while preserving solution quality in large-scale power system optimization under strict time limits.
\subsection{Contribution}
Building on our prior research \cite{TPEC, IEEE-TIA}, this paper presents a GPU-accelerated solver for the unit commitment (UC) problem in large-scale power systems. The solver leverages the Primal–Dual Hybrid Gradient (PDHG) algorithm to efficiently solve the relaxed linear subproblem, improving bound estimation as well as crossover and branch-and-bound convergence. By exploiting GPU acceleration, the approach offers a scalable and practical framework for real-time large-scale UC applications. The method is validated on 4224-, 6049-, and 6717-bus networks, demonstrating substantial computational speed-ups while preserving solution quality. 




\begin{figure*}
\centerline{\includegraphics[scale=0.15]{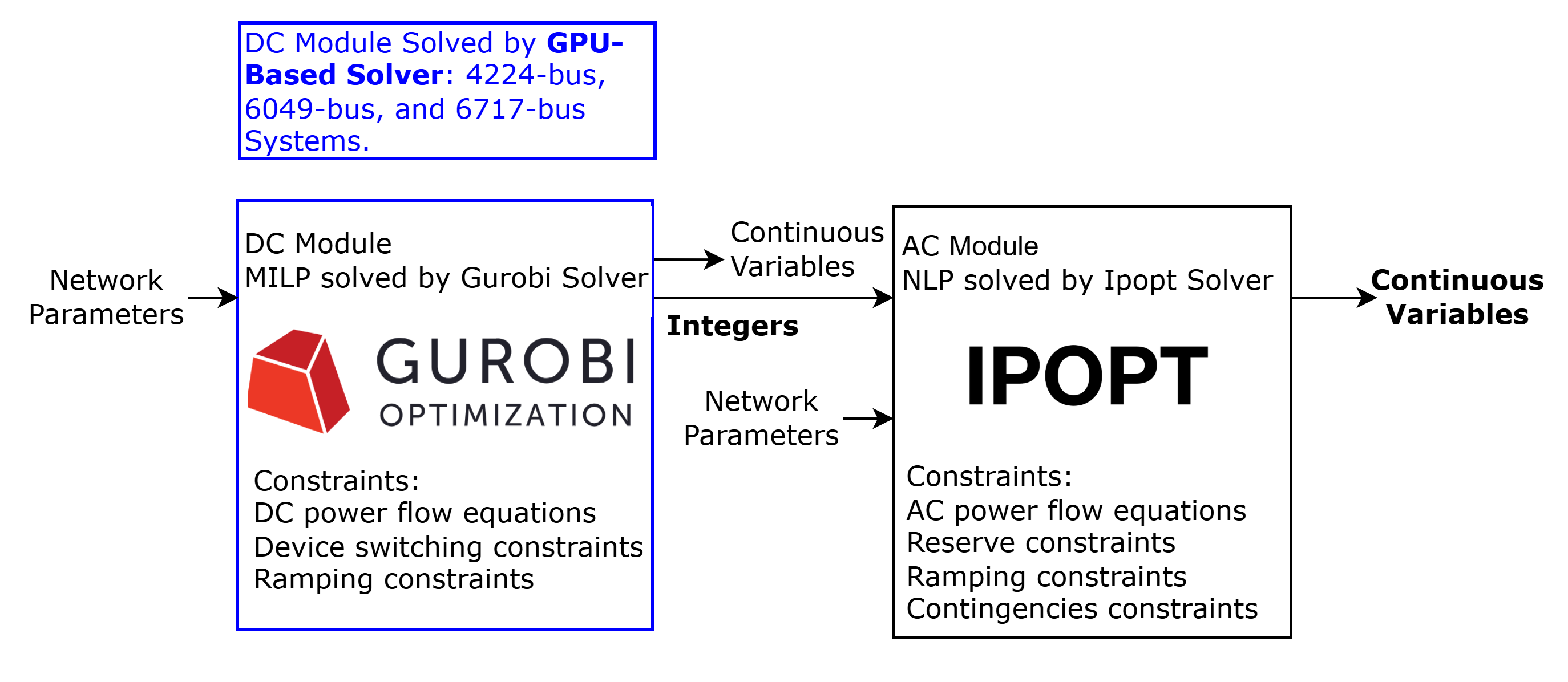}}
\caption{Our two-stage solution framework integrates the DC and AC modules presented in \cite{TPEC, IEEE-TIA}. In this work, we improve the computational speed of the DC module (first stage) using GPU-based optimization. The GPU-based solver accelerates computations for large-scale networks, including the 4224-, 6049-, and 6717-bus systems.}
\vspace{-.5cm}
\label{fig0}
\end{figure*}

\section{Problem Setup}

\subsection{Problem Formulation}
This study adopts the problem formulation provided by the DOE Advanced Research Projects Agency-Energy (ARPA-E) for the GO challenge 3. The problem is highly detailed, with the primary formulation spanning 53 pages \cite{GO3}, accompanied by an additional 22 pages of supplementary material. It involves a comprehensive set of optimization variables and constraints, such as limits on real and reactive power at individual bus nodes, voltage bounds, zone-specific reserve requirements, device operational statuses, switching complexities, considerations for device downtime and uptime, the total number of device starts over multiple intervals, ramping constraints, and minimum/maximum energy limits across time periods 
(refer to Fig. 1 of \cite{TPEC}). 
Notably, all constraints are applied under both base and contingency scenarios. The contingency conditions ensure the system maintains reliability under unexpected events and outages.

\subsection{Dataset}
This study utilizes the ARPA-E dataset provided in \cite{go_data}. The 4224-bus, 6049-bus, and 6717-bus networks from Trial 3 were adapted for our analysis. These networks correspond to different scheduling horizons: 1 day ahead (Division 1), 2 days ahead (Division 2), and 1 week ahead (Division 3). The optimization model for the 4224-bus network with Division 1 includes approximately 1 million continuous decision variables, around 116,000 binary variables, and roughly 1 million constraints. The complete counts of variables and constraints for all networks are summarized in Table \ref{tab01}. The values in Table \ref{tab01}  correspond to Division 1 (18 steps). To obtain Division 2 values (48 steps), divide by 18 and multiply by 48. For Division 3 (42 steps), divide by 18 and multiply by 42. In total, 57 scenarios across the three networks and three scheduling horizons were tested in this study.

\begin{table}[b]
\caption{Problem scale for different network sizes. The values correspond to Division 1 (18 steps). To obtain Division 2 values (48 steps), divide by 18 and multiply by 48. For Division 3 (42 steps), divide by 18 and multiply by 42.}
\vspace{-.25cm}
\begin{center}
\begin{tabular}{|l|c|c|c|}
\hline
\textbf{Parameter} & \multicolumn{3}{c|}{\textbf{Network}} \\
\hline
 & \textbf{4224-Bus} & \textbf{6049-Bus} & \textbf{6717-Bus} \\
\hline
\# Binary Variables & 116,154 & 67,932 & 104,868 \\
\hline
\# Continuous Variables & 966,220 & 1,673,301 & 2,165,998 \\
\hline
\# Constraints & 1,168,114 & 1,580,545  & 2,003,974 \\
\hline
\end{tabular}
\label{tab01}
\end{center}
\end{table}

\subsection{Reference Model}
Our previous model introduced in \cite{TPEC, IEEE-TIA} demonstrated excellent performance, with an average scaled score above 0.98 for a range of network sizes and control settings. The scaled score reflects how closely the obtained solution approaches the best-known outcome, and consistently achieving 0.98 over roughly 300 scenarios highlights the reliability of the prior approach.

In this work, we improve the speed of convergence of the unit commitment module by leveraging GPU-based optimization, solved using the Gurobi solver, as depicted in Figure \ref{fig0}.

\subsection{Hardware and Computing Setup}
Running the PDHG algorithm on a GPU with the Gurobi solver requires NVIDIA hardware (preferably an NVIDIA H100 GPU, as noted in the Gurobi documentation \cite{gurobi}), CUDA version 12.9, and a Linux 64-bit operating system. The GPU-based solver utilizes NVIDIA’s cuOpt engine. While Gurobi is typically supported and widely used on Windows platforms, GPU acceleration is currently available only on Linux systems.

The cuOpt engine is NVIDIA’s GPU-accelerated optimization library designed for large-scale linear and mixed-integer programming. It leverages highly parallel numerical kernels to execute optimization steps efficiently on GPUs, providing substantial speed-ups compared to traditional CPU solvers. In Gurobi’s GPU mode, cuOpt handles the PDHG iterations, while Gurobi oversees higher-level tasks such as presolve, crossover, and branch-and-bound.

In this study, we employed an NVIDIA A800 (40 GB) active GPU configured under Linux64 through the Windows Subsystem for Linux (WSL). The Gurobi 13.0.0 beta 2 version was used. 
The computing node was equipped with an Intel® Xeon® W9-3495X processor with 56 cores (112 threads), operating at a base frequency of 1.90 GHz and a maximum turbo boost of 4.80 GHz, along with 256 GB of RAM to handle the size and complexity of the datasets. Gurobi was configured to use 64 threads (the default is 32). Further increasing the thread count to the maximum available did not improve performance, likely due to memory bandwidth constraints and synchronization overhead.

Table~\ref{tab:hardware} summarizes the hardware and software specifications used in our experiments. The GPU and multicore CPU operated collaboratively to accelerate the optimization process. The PDHG algorithm was executed on the GPU, while the remaining steps—such as crossover and branch-and-bound—benefited from multithreading on the CPU. In contrast, when linear programming and barrier methods are employed instead of PDHG, they primarily benefit from CPU-based multithreading rather than GPU acceleration.

\begin{table}[h]
\centering
\caption{Hardware and Software Specifications}
\label{tab:hardware}
\begin{tabular}{|l|l|}
\hline
\textbf{Component} & \textbf{Specification} \\ \hline
GPU & NVIDIA A800, 40 GB VRAM \\ \hline
GPU Engine & NVIDIA cuOpt library \\ \hline
CPU & Intel® Xeon® W9-3495X, 56 cores\\ 
 & / 112 threads, 1.90–4.80 GHz \\ \hline
RAM & 256 GB \\ \hline
Operating System & Linux 64-bit (via WSL on Windows) \\ \hline
Python Version & 3.11 \\ \hline
Gurobi Version & 13.0.0 beta 2\\ \hline
CUDA Version & 12.9 \\ \hline
Gurobi Threads & 64 (default 32) \\ \hline
\end{tabular}
\end{table}







\subsection{PDHG Background}
For many years, the field of linear programming (LP) has been dominated by two primary algorithms: the Simplex Method and Interior Point Methods (IPMs). These approaches have been extensively refined and remain powerful and reliable for a wide range of optimization problems. However, adapting them for GPU acceleration is challenging due to their dependence on complex matrix operations and factorization routines.

The Primal-Dual Hybrid Gradient (PDHG) method, also known as PDLP \cite{app2021} when adapted for LPs, offers a first-order alternative designed for large-scale problems. Unlike traditional methods, PDHG relies on gradient-based updates and lightweight operations such as sparse matrix-vector multiplications, which are well-suited for GPU implementation. The algorithm builds on the primal–dual formulation of LPs, iteratively updating both primal and dual variables to solve saddle-point problems. 
This framework was initially introduced by Chambolle and Pock \cite{Chambolle2011}, and subsequent enhancements by Applegate et al. \cite{app2021} incorporated adaptive step size control and restart mechanisms to accelerate convergence.

\noindent  A standard LP and its dual can be expressed as:

\[
\text{Primal: } \min_x \; c^\top x \quad \text{s.t. } Ax \le b, \; x \ge 0,
\]  
\[
\text{Dual: } \max_y \; b^\top y \quad \text{s.t. } A^\top y \ge c, \; y \ge 0.
\]

This pair of problems can be equivalently written in a saddle-point form as:

\[
\min_x \max_y \; L(x,y) = c^\top x - y^\top A x + b^\top y.
\]

The PDHG method solves this by performing iterative updates:

\[
x_{k+1} = \mathrm{proj}_{\mathbb{R}^n_+}\Big(x_k + \tau A^\top y_k - \tau c \Big),
\]  
\[
y_{k+1} = y_k - \sigma A \big(2 x_{k+1} - x_k\big) + \sigma b,
\]

where $\tau$ and $\sigma$ denote the primal and dual step sizes, respectively. Building on this foundation, Applegate et al. \cite{app2021} introduced practical modifications—including adaptive step size selection and restart strategies—that enhance the performance of PDHG on large LPs.

\section{Results and Discussion}

\subsection{Evaluation of Solver Runtime Performance}

\begin{figure}
    \centering
    \includegraphics[width=0.9\linewidth]{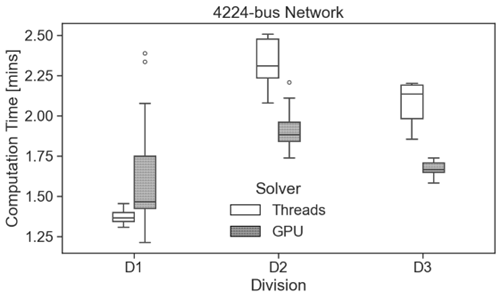}
    \caption{Computation time distribution for the 4224-bus network (Divisions 1–3).  The GPU-based solver achieves faster convergence and lower variance in Divisions 2 and 3.}
    \label{fig:4224}
\end{figure}

\begin{figure}
    \centering
    \includegraphics[width=0.9\linewidth]{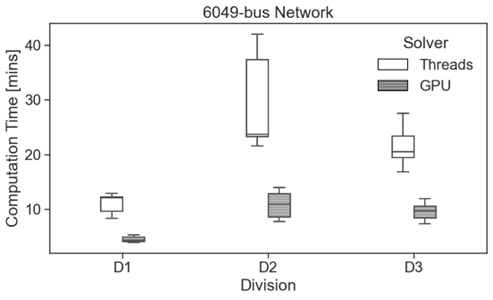}
    \caption{Computation time distribution for the 6049-bus network (Divisions 1–3). The GPU-based solver achieves faster convergence and lower variance across all Divisions.}
    \label{fig:6049}
\end{figure}

\begin{figure}
    \centering
    \includegraphics[width=0.9\linewidth]{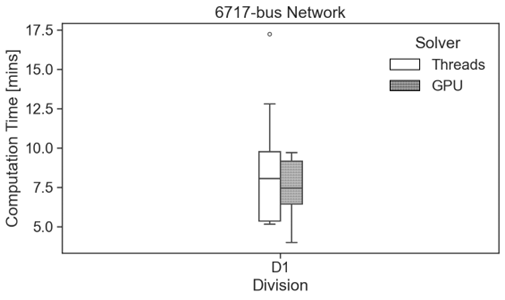}
    \caption{Computation time distribution for the 6717-bus network (Division 1). The GPU-based solver shows improved convergence.}
    \label{fig:6717}
\end{figure}

The computation times for the 4224-, 6049-, and 6717-bus networks were evaluated across multiple scenarios, Divisions 1–3, and solver configurations, comparing GPU-based PDHG against multi-threaded CPU solvers using the barrier method. Figures \ref{fig:4224}, \ref{fig:6049}, and \ref{fig:6717} show boxplots of the total runtime for each network and division. Overall, the GPU solver provides faster computation, particularly for larger networks and computationally intensive divisions (e.g., long control horizons in Divisions 2 and 3).

In the 4224-bus network, GPU-based PDHG outperforms CPU threads in Divisions D2 and D3. The most significant speedups are observed in the 6049-bus network, where CPU runtimes for certain divisions spike, while PDHG consistently reduces total computation time and exhibits lower variance. For the 6717-bus network, GPU-based PDHG also demonstrates overall superior performance compared to CPU threads.

These improvements are largely attributed to the higher-quality solutions produced by PDHG, which provide a superior warm start for the crossover and branch-and-bound stages, thereby accelerating convergence. While CPU-based barrier solvers benefit from multi-threading, they cannot exploit the parallelism offered by GPU execution. The boxplots also reveal scenario-dependent variability, reflecting differences in convergence behavior between PDHG and the barrier method. Overall, GPU-based PDHG shows the greatest advantage for large-scale and computationally demanding problems, such as Divisions 2 and 3 of the 6049-bus system, providing faster and more consistent runtimes with reduced variability.

\vspace{-0.2cm}

\subsection{PDHG Execution \& Crossover Refinement}

The PDHG iteration log is similar to that of the barrier method, displaying the current primal and dual objective values, primal and dual residuals, and complementarity. As the iterations progress, the gap between the primal and dual objectives narrows, and both the residuals and complementarity decrease.

The termination criteria used in this study were set to their default values (see Termination for more details in \cite{gurobi}). Once the criteria are satisfied, the crossover step converts the PDHG solution into a basic feasible solution.

Unlike the Simplex method, which directly produces a basic feasible solution, PDHG and interior algorithms generate dense solutions that may slightly violate feasibility or optimality conditions. Rather than switching to the Simplex algorithm—which tends to be slower and less suitable for GPU acceleration—the crossover phase efficiently refines the output of PDHG or interior algorithms into a clean, sparse basic solution. It should be noted that only the PDHG algorithm itself is executed on the GPU.

It is instructive to experiment with different termination tolerances. For example, increasing the PDHG termination tolerances, such as the relative LP duality gap, reduces the computational effort spent within the PDHG iterations. However, this comes at the cost of solution accuracy, requiring the crossover algorithm to spend more time refining the approximate solution into a high-quality feasible one.

By design, first-order methods such as PDHG converge linearly, meaning they rarely achieve highly accurate solutions at absolute tolerances. Therefore, PDHG is typically terminated once relative constraint violations fall below a practical threshold. To improve accuracy, Gurobi automatically invokes the crossover algorithm, which refines PDHG’s approximate solution into a high-accuracy one. Although crossover and branch-and-bound still run on the CPU, they benefit from the improved solution quality produced by the GPU-based PDHG. Specifically, PDHG provides a better warm start for crossover, a higher-quality initial basis, and tighter LP relaxation bounds for MILP problems. Consequently, even CPU-only stages complete faster, as they begin from a more accurate and well-conditioned solution.
\begin{figure*}[ht]
    \centering
    \begin{subfigure}[b]{0.48\linewidth}
        \centering
        \includegraphics[width=0.8\linewidth]{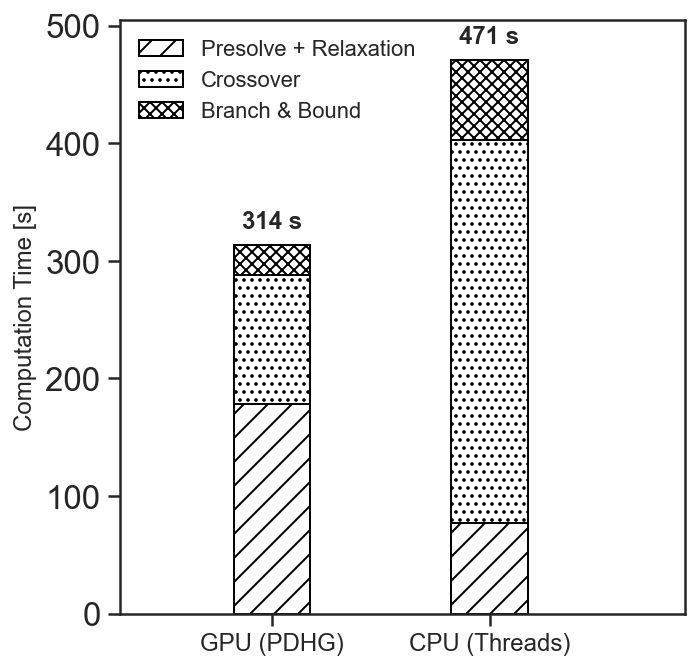}
        \caption{Random Scenario 1}
        \label{fig:Time_breakdown1}
    \end{subfigure}
    \hfill
    \begin{subfigure}[b]{0.48\linewidth}
        \centering
        \includegraphics[width=0.8\linewidth]{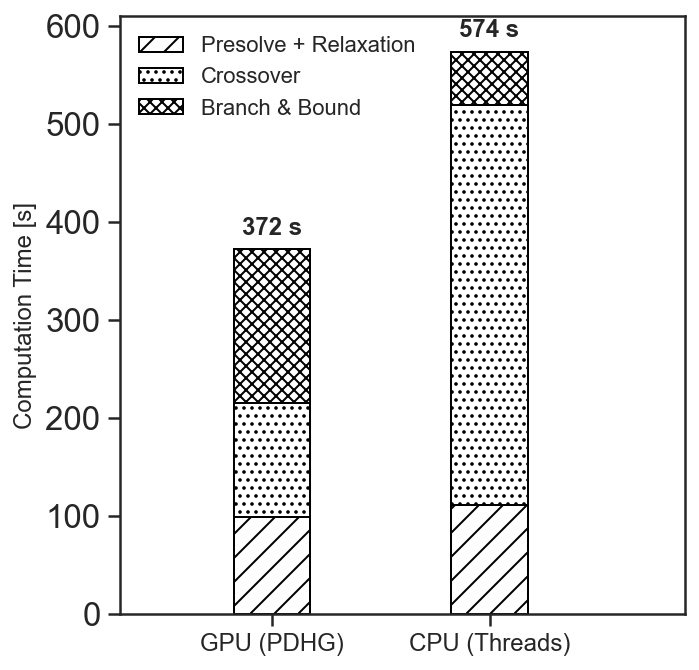}
        \caption{Random Scenario 2}
        \label{fig:Time_breakdown2}
    \end{subfigure}
    \caption{Comparison of computation times between GPU-based PDHG and multi-threaded CPU solvers for the 6049-bus system across two random scenarios. Each bar represents the total runtime, decomposed into three stages: (1) presolve and relaxation, (2) crossover, and (3) branch-and-bound. For the CPU solver, the relaxation stage corresponds to the barrier method, while for the GPU solver, it corresponds to PDHG. The GPU-based relaxation consistently provides a higher-quality initial solution, reducing overall computation time by accelerating convergence in the crossover and/or branch-and-bound stages.}
    \label{fig:Time_breakdown_combined}
\end{figure*}

Figure \ref{fig:Time_breakdown_combined} illustrates the computation time breakdown for solving the problem using GPU-based PDHG and multi-threaded CPU solvers for two random scenarios of the 6049-bus system. In both cases, the GPU solver demonstrates a clear advantage in overall runtime, particularly in the crossover stage, which dominates CPU computation time. In the first scenario, the presolve and relaxation phases of the CPU solver—corresponding to the barrier method—are faster on the CPU, but the GPU-based relaxation using PDHG significantly reduces the crossover and branch-and-bound stages, resulting in a lower total runtime. In the second scenario, the GPU is slightly faster in the PDHG relaxation phase, and although the branch-and-bound stage is slower, the substantial reduction in crossover time still enables the GPU-based approach to achieve a lower total computation time. This improvement is attributed to the high-quality solution obtained by PDHG, which provides a better warm start and a more accurate initial basis for the crossover, accelerating this stage. 

Specifically, the crossover stage starting from the GPU-based relaxation is more than $2.96\times$ faster than when starting from the CPU-based relaxation for these two scenarios. Overall, the GPU-based solver is approximately $1.5\times$ faster than the CPU-based solver across both scenarios.

\subsection{Speedup Summary}
Solver runtimes highlight the computational benefits of utilizing GPUs to solve the relaxed MILP for the unit commitment problem. Tables~\ref{tab:max_times} and~\ref{tab:avg_speedup} summarize the maximum solution times and the speedups achieved with the GPU-based solver. The GPU solver reduces the maximum solution time and overall computation time for the 4224-bus system in divisions 2 and 3 (speedups $>1.2\times$), for the 6049-bus system across all divisions (speedups $>2.27\times$), and for the 6717-bus system (speedup $=1.23\times$). Notably, for the 6049-bus system in division 2, the GPU solver decreases the maximum solution time from 42 minutes to 14 minutes.

\subsection{Solution Quality}
The GPU-accelerated solver maintains the quality of the solutions, as  summarized in Table \ref{tab:obj_score_distribution}. This table presents the distribution of scores, defined as the ratio between the GPU-computed objective and the thread-based objective. Since the objective corresponds to the maximized market surplus \cite{GO3}, a score greater than 1 indicates an improved solution.

Out of 57 test scenarios, 36 scenarios achieved a score of 1 or higher, while the remaining 21 scenarios had scores slightly below 1, yet still above 0.999. The overall mean score across all test scenarios is approximately 0.9999, demonstrating that the GPU-based solver effectively preserves high-quality solutions.

 \subsection{Mitigating Threads Spin/Wait Overhead}

In the CPU-based runs, the solver may execute the primal simplex, dual simplex, and barrier algorithms in parallel using multithreading. However, once one of these algorithms converges, Gurobi continues to wait for the others to finish, which can introduce spin/wait overhead and result in idle CPU cycles. This behavior was observed in the 6049-bus system, as reported in Table~\ref{tab:max_times}. To mitigate this effect and more accurately reflect the computational performance of the CPU-based solver, we enabled the setting that terminates the remaining algorithms once a converged solution is reached, thereby reducing the spin/wait behavior. This adjustment improved the Threads-based runtimes and ensured a fair and consistent comparison with the GPU-based solver. The speedup achieved by the GPU over the CPU-based solver after addressing the spin/wait time is shown in Table~\ref{tab:avg_speedup_updated}, with the average speedup ranging from 1.31× to 1.93×. Overall, the GPU-based solver continues to outperform the CPU-based solver.

The improvement in average computation time when enabling the setting to terminate remaining algorithms once one completes is summarized in Table~\ref{tab:CPU_avg_speedup}. This indicates that the adjustment leads to more efficient CPU utilization. These results highlight that careful management of multithreaded solver behavior can substantially improve performance and provide a more accurate and fair comparison with GPU-based implementations.

\subsection{Future Outlook}
Given the rapid evolution of GPU hardware and supporting libraries—such as NVIDIA’s cuDSS and cuSparse, which were utilized in this GPU-based implementations—the performance outcomes reported here are expected to improve over time, as the solver stands to benefit directly from advances in parallel sparse linear algebra kernels and memory throughput.

\begin{table}[h!]
\caption{\textcolor{black}{Maximum solution times for Threads-based and GPU-based solvers. Arrows indicate whether GPU is faster ($\downarrow$, blue) or slower ($\uparrow$, red) than Threads. Time in minutes.}}
\centering
\small
\setlength{\tabcolsep}{4pt}
\begin{tabular}{|l|l|c|c|}
\hline
\textbf{Network} & \textbf{Division (D)} & \multicolumn{2}{c|}{\textbf{Maximum Time [mins]}}  \\
\cline{3-4}
 & & Threads & GPU  \\
\hline
4224-bus & D1 & 1.46 & 2.39 \textcolor{red}{$\uparrow$} \\
     & D2 & 2.51 & 2.21 \textcolor{blue}{$\downarrow$} \\
     & D3 & 2.20 & 1.74 \textcolor{blue}{$\downarrow$} \\
\hline
6049-bus & D1 & 12.94 & 5.41 \textcolor{blue}{$\downarrow$} \\
     & D2 & 42.03 & 14.0 \textcolor{blue}{$\downarrow$} \\
     & D3 & 27.57 & 12.0 \textcolor{blue}{$\downarrow$} \\
\hline
6717-bus & D1 & 17.24 & 9.71 \textcolor{blue}{$\downarrow$} \\
\hline
\end{tabular}
\label{tab:max_times}
\end{table}

\begin{table}[h!]
\caption{\textcolor{black}{Average computation times for Threads-based and GPU-based solvers and resulting speed-up. Colors indicate whether GPU is faster (blue) or slower (red) than Threads. Time in minutes.}}
\centering
\small
\setlength{\tabcolsep}{4pt}
\begin{tabular}{|l|l|c|c|c|}
\hline
\textbf{Network} & \textbf{Division (D)} & \multicolumn{2}{c|}{\textbf{Avg. Time [mins]}} & \textbf{Avg.} \\
\cline{3-4}
 & & Threads & GPU&\textbf{Speed-up} \\
\hline
4224-bus &  D1 & 1.37 & 1.61 & \textcolor{red}{$0.85\texttt{x}$} \\
     & D2 & 2.34 & 1.92 & \textcolor{blue}{$1.22\texttt{x}$} \\
     & D3 & 2.08 & 1.67 & \textcolor{blue}{$1.25\texttt{x}$} \\
\hline
6049-bus & D1 & 11.13 & 4.56 & \textcolor{blue}{$2.44\texttt{x}$} \\
     & D2 & 29.36 & 10.19 & \textcolor{blue}{$2.88\texttt{x}$} \\
     & D3 & 21.85 & 9.64 & \textcolor{blue}{$2.27\texttt{x}$} \\
\hline
6717-bus & D1 & 8.76 & 7.12 & \textcolor{blue}{$1.23\texttt{x}$} \\
\hline
\end{tabular}
\label{tab:avg_speedup}
\end{table}

\begin{table}[ht]
\caption{Distribution of scores for GPU objective. Score $=$ Objective (GPU) $/$ Objective (Threads). Out of 57 scenarios, 36 scenarios achieve a score of at least 1, while 21 scenarios have scores slightly below 1 but above 0.999, indicating that the GPU solver maintains solution quality.}
\centering
\begin{tabular}{|l|c|}
\hline
\textbf{Score Range} & \textbf{Number of Scenarios} \\
\hline
$0.999 <$ Score $< 1$ & \textcolor{red}{21} \\
\hline
\textbf{Score $\geq 1$} & \textcolor{blue}{36} \\
\hline
\hline
\textbf{Average Score} & $\approx$ \textcolor{red}{$0.9999$} \\
\hline
\end{tabular}
\label{tab:obj_score_distribution}
\end{table}

\begin{table}[h!]
\caption{\textcolor{black}{Average computation times for Threads-based and GPU-based solvers and resulting speed-up, for the 6049-bus system. The Threads-based solver runtime was improved by addressing CPU spin/wait effects. Colors indicate whether GPU is faster (blue) or slower (red) than Threads.
}}
\centering
\small
\setlength{\tabcolsep}{4pt}
\begin{tabular}{|l|l|c|c|c|}
\hline
\textbf{Division (D)} & \multicolumn{2}{c|}{\textbf{Avg. Time [mins]}} & \textbf{Avg.} \\
\cline{2-3}
 & Threads & GPU&\textbf{Speed-up} \\
\hline

 D1 & 5.97 & 4.56 & \textcolor{blue}{$1.31\texttt{x}$} \\
 D2 & 19.7 & 10.19 & \textcolor{blue}{$1.93\texttt{x}$} \\
D3 & 16.06 & 9.64 & \textcolor{blue}{$1.67\texttt{x}$} \\
\hline
\end{tabular}
\label{tab:avg_speedup_updated}
\end{table}

\begin{table}[h!]
\caption{\textcolor{black}{Speedups achieved for the threads-based solver runtime by addressing CPU spin/wait effects for the 6049-bus system.}}
\centering
\small
\setlength{\tabcolsep}{4pt}
\begin{tabular}{|l|c|}
\hline
\textbf{Division (D)} & \textbf{Avg. Speed-up} \\
\hline

 D1  & \textcolor{blue}{$1.86\texttt{x}$} \\
 D2  & \textcolor{blue}{$1.49\texttt{x}$} \\
D3 & \textcolor{blue}{$1.36\texttt{x}$} \\
\hline
\end{tabular}
\vspace{-.3cm}
\label{tab:CPU_avg_speedup}
\end{table}

\section{Conclusion}
\vspace{-.15cm}
This paper presents a GPU-accelerated approach for solving the unit commitment problem in large-scale power systems. By leveraging the Primal–Dual Hybrid Gradient (PDHG) algorithm on GPUs, the proposed solver achieves faster convergence in the relaxed linear subproblem, providing high-quality warm starts for crossover and branch-and-bound procedures. Validation on 4224-, 6049-, and 6717-bus networks demonstrates significant reductions in computation time—up to $2.88\times$ speed-up—while preserving solution quality, with mean scores near 0.9999 relative to CPU-based methods. These results highlight the potential of GPU-based optimization to enhance the efficiency and scalability of large-scale mixed-integer linear programming for power system operations, enabling faster and more reliable unit commitment in modern grids. As GPU architectures and sparse optimization libraries continue to advance, the performance and scalability of the proposed solver are expected to further improve.

\section*{Acknowledgment}
\vspace{-.2cm}
Artificial intelligence tools were used to enhance the clarity and readability of the manuscript.
We want to thank the support of US ARPA-E (\#DE-AR0001646), and NSF (\#2442420 and \#2313768).

\end{document}